\documentclass[preprint,12pt]{elsarticle}
\usepackage{amsmath,amssymb}
\usepackage{mathptmx}
\usepackage{epsfig,latexsym}
\usepackage{graphics}
\usepackage[latin2]{inputenc}

\newcommand{\R}{\mathbb{R}}

\newcommand{\RN}{\mathbb{R}^N}
\newcommand{\RM}{\mathbb{R}^M}

\newcommand{\ub}{\mathbf{u}}

\newcommand{\xb}{\mathbf{x}}

\newcommand{\PDE}{partial differential equation}

\renewcommand{\emph}{\textbf}

\newcommand{\torol}[1]{}

\newtheorem{rmk}{Remark}
\newtheorem{thm}{Theorem}
\newtheorem{prp}{Proposition}
\newcommand{\Pf}[1]{\textit{Proof. }#1\(\square\)\par}

\journal{Computational and Mathematical Modeling}

\begin{document}
\begin{frontmatter}

\title{Application of Operator Splitting to Solve Reaction--Diffusion Equations}
\author{Tam\'as Ladics}
\address{Szent Istv\'{a}n University, Ybl Mikl\'{o}s College of Building \fnref{bme}}
\fntext[bme]{Research done at the Department of Mathematical Analysis,
Budapest University of Technology and Economics,
Budapest, H-1111 Egry J. u. 1., HUNGARY}

\begin{abstract}
Approximate solutions of the Fisher equation obtained by different splitting methods are investigated. The error of this nonlinear problem is analyzed. The order of different splitting methods coupled with numerical methods of different order is calculated numerically and symbolically.
\end{abstract}


\end{frontmatter}

\section{Introduction}
Splitting methods have been fruitfully used to solve large systems of partial differential equations. To find the exact solution of a given problem in practice is usually impossible. We can use numerical methods to obtain an approximate solution of the equations, although the discretized model can be still very difficult to solve. Reaction-diffusion models or models of transport processes have a structure that allows a natural decomposition of the equations, thus provide the opportunity to apply operator splitting schemes. Splitting methods help us reduce the complexity of the system and reduce computational time. With splitting it is possible to handle stiff terms separately and to solve each subproblem with a suitable numerical method chosen to the corresponding operator. To solve a problem in practice we use operator splitting and numerical schemes which we will call the \textit{combined method}. 
The use of operator splitting as well as the numerical methods result in some error in the solution. The error generated purely by splitting is called \textit{splitting error}. This is the difference of the exact solution and the approximate solution obtained by splitting (assumed that we know the exact solutions of the subproblems). Combined methods can generate both splitting error and \textit{numerical error}. The study of this common effect on the solution is our main concern in this paper. Detailed study on the interaction of operator splitting and numerical schemes for \textit{linear} problems can be found in \cite{csomosfarago}. They classify the errors that can occur using splitting methods and numerical schemes, give theoretical and numerical results on the order of the combined method for linear problems. Our aim is to characterize the error of this combined method therefore we calculate the order of the combined method for a \textit{nonlinear} problem. We analyze the order of the error in the light of the characteristics of the splitting error and the \textit{numerical error}.

\cite{sanz} and \cite{lanserverwer} discuss the splitting error in a general framework.
The effect of operator splitting on the wave solutions of the Fisher equation is investigated by \cite{simpson}. 
Our aim here is to rigorously analyze the interaction of splitting error and numerical error in the case of a nonlinear problem: the Fisher equation. The structure of our paper is as follows. In section 2 we introduce the basic idea of operator splitting in a general frame. In Section 3 we introduce the Fisher equation and recall some known results on it. We show how we apply splitting to solve the Fisher equation. In Section 4 the splitting error is analyzed for reaction-diffusion problems in full generality. In Section 5 we calculate the order of the combined method for nonlinear problems in general. Section 6 contains the numerical results on the Fisher equation.

\section{Operator splitting} 

Let us consider the following abstract Cauchy problem:
\begin{equation}\label{os}
U^\prime(t)=A(U(t)) \quad U(0)=U_0
\end{equation}
with $U_0\in X$ arbitrary. The set $X$ is usually a space of functions with certain properties, $U(t)\in X$ for every $t\geqslant 0$ and $A:X\to X.$
Suppose that $A$ can be written as the sum of two operators: $A=A_1+A_2.$

The most simple type of operator splitting is the sequential splitting. In this case the split problem is:
\begin{equation}\label{os1}
U_1^\prime(t)=A_1(U_1(t)) \quad U_1(0)=U_0
\end{equation}

\begin{equation}\label{os2}
U_2^\prime(t)=A_2(U_2(t)) \quad
U_2(0)=U_1(\tau).
\end{equation}

The basic idea of splitting is to decompose the operator on the
right hand side into the sum of simpler operators, and to solve the
subproblems corresponding to the operators successively in each time
step. More precisely, we solve the equation only with operator
$A_1$ until time $\tau$ (as if only the subprocess represented by
$A_1$ were present) and the solution in time $\tau$ will be the
initial condition of the equation with $A_2$. It means that we
return to the initial time and solve the equation with $A_2$ as
well. The solution of the second equation in time $\tau$ is called the
approximate solution of the original problem in time $\tau.$ This
procedure is then repeated on the interval $[\tau, 2\tau]$ etc.
Thus, the simpler subproblems are connected to each other through
the initial conditions. It is clear, that the numerical treatment of
the separate subproblems is simpler. The most significant advantage
of splitting is that we can exploit the special properties of the
operators of the different subproblems and apply the most suitable
numerical method for each of them. Thus we can obtain a more precise
solution in a shorter time.
\\
We remark that the method can be used fruitfully in large models, for
example global models of air pollution transport, or combustion or metabolic models, where the number
of predicted variables is large and the number of the processes
represented in the models is large. We refer to three works on air pollution models with application of operator splitting of \cite{lagzi1}.

\subsection{Splitting schemes}

We can define the different splitting methods by solving the subproblems successively in different orders and for different time lengths. The above described simplest scheme is called  \textit{sequential} splitting (SEQ). We solve the subproblems one after another using the same time length $\tau,$ schematically $S_2(\tau) S_1(\tau),$ where $S_i$ is the corresponding solution operator. In \textit{Marchuk--Strang} splitting (MS) we usually solve the subproblem with $A_1$ with time length $\tau/2$ then solve the other one with $A_2$ for time length
$\tau$ and solve with $A_1$ again for time length $\tau/2$. Schematically $S_1(\tau/2) S_2(\tau) S_1(\tau/2).$ We usually chose $A_2$ to be the operator representing chemical reactions. In general the operator that is stiff or nonlinear. In this given order we only need to solve the second subproblem once which can be of importance given the operator's properties. In \textit{weighted sequential} splitting the solution in the next time step is a weighted average of the results of the two possible sequential splittings $S_1(\tau) S_2(\tau)$ and $S_2(\tau) S_1(\tau)$. In the special case of \textit{symmetrically weighted} splitting (SW) we take the arithmetic mean of the results: $(S_1(\tau) S_2(\tau)+S_2(\tau) S_1(\tau))/2.$ The extra work with MS and SW splittings benefits in second order accuracy compared to the first order of SEQ splitting. The nonsymmetric weighted splitting is of order one. In later sections we investigate the SEQ, the MS and the SW splittings coupled with four different numerical methods, all of different orders.

\subsection{Splitting error, order of splitting}
We perform a semidiscretization on (\ref{os}) in an equidistant manner with time step $\tau.$ If we know the exact solutions of (\ref{os1}) and (\ref{os2}) we can generate an approximate solution to the original full problem (\ref{os}) in which error originated only from operator splitting can arise. If we denote the exact solution by $U$ and the approximate solution by $\tilde{U}$ then the local error of operator splitting is 
\[
E(\tau):=U(\tau)-\tilde{U}(\tau)
\]
Both solutions start from the common initial value and after time $\tau$ the difference $E(\tau)$ is called \textit{splitting error}.
Naturally the splitting error can be defined at any point of time during integration, if $U(t)=\tilde{U}(t)$ then $U(t+\tau)-\tilde{U}(t+\tau)$ is the local error at $t.$

For linear operators it is easy to show by Taylor expansion that the local order of SEQ equals $1$ since the error becomes $E(\tau)=K \tau^2+O(\tau^3).$ For nonlinear operators we need the definition of the \textit{Lie-operator} and we can perform  the analysis with Taylor expansion using the Lie-operators. We refer to \cite{lanserverwer} for detailed derivation of the nonlinear case. From the literature on operator splitting it is well known that the MS provides second order accuracy, so does the SW splitting, see \cite{faragohavasi}.

\subsection{Splitting of reaction--diffusion equations}
In the case of reaction--diffusion equations there is a natural decomposition of (\ref{os}). The operator $A_1$ represents the process of diffusion and $A_2$ the chemical reactions. In this case $A_1$ is a linear and unbounded operator, $A_2$ is usually a nonlinear operator. If we have $M$ species and $N$ denotes the spatial dimension of the problem then $U\in C^1(\mathbb{R}^+, C^2(\RN,\RM))$ that is for a given time point $t\in\R^+$ the function $U$ maps the concentration of all species for every given point in space $\R^N,$ so it is a function of $\xb\in \R^N.$ In other words $U:t\mapsto (u_1(t,\xb),...,u_M(t,\xb)),$ where $u_i(t,\xb)$ is the concentration of the $i$th species that is the spatial distribution of the $i$th species. Then the function space $X:=C^1(\mathbb{R}^+,C^2(\RN,\RM)).$

\section{The Fisher equation}

The Fisher equation is:
\begin{equation}\label{fisher}
\left\{
\begin{array}{lll}
\partial_tu (t,x)& = & \partial^2_xu(t,x)+u(t,x)(1-u(t,x)) \qquad x\in \R, t> 0\\
u(0,x) & = & \eta(x). \\
\end{array}
\right.
\end{equation}

There is only one chemical species present and one spatial variable here. This equation was originally derived to describe the propagation of a gene in a population \cite{fisher}. It is one of the simplest nonlinear models for reaction-diffusion equations. Such equations occur, e.g., in combustion, mass transfer, crystallization, plasma physics, and in general phase transition problems. See a discussion on reaction-diffusion models in \cite{erditoth} and \cite{murray}.
For the initial condition:
\[
u(0,x)=\frac{1}{(1+k\,\mathrm{exp}(x/\sqrt{6}))^2}
\]
wave form solution of the equation is known:
\[
u(t,x)=\frac{1}{\left(1+k\,\mathrm{exp(-\frac{5}{6}t+\frac{1}{6}\sqrt{6}x)}\right)^2}
\]
and
for:
\[
u(0,x)=\frac{1}{(1+k\,\mathrm{exp}(-x/\sqrt{6}))^2},
\]
\[
u(t,x)=\frac{1}{(1+k\,\mathrm{exp(-\frac{5}{6}t-\frac{1}{6}\sqrt{6}x)})^2}
\]
We investigate three different splitting methods applied in the solution of this equation. A
natural way to split the Fisher equation is to decompose it into two
subproblems: one for the diffusion and one that corresponds to the
reaction part of the right hand side. Thus the definitions of the
subproblems are:
\begin{equation}\label{subpra}
\left\{{\setlength\arraycolsep{.13889em}
\begin{array}{lll}
\partial_t u_1 (t,x) & = & \partial^2_x u_1 (t,x)\\
u_1(0,x) & = & \eta_1(x) \\
\end{array}}
\right.
\end{equation}
\begin{equation}\label{subprb}
\left\{
\begin{array}{lll}
\partial_tu_2 (t,x)& = & u_2(t,x)(1-u_2(t,x)) \\
u_2(0,x) & = & \eta_2(x), \\
\end{array}
\right.
\end{equation}
where the initial condition $\eta_2(x)=u_1(\tau,x)$ connects the equations.
The lower indexes help distinguish the solutions of the different
problems. We define the operators $A_1$ and $A_2$ as follows: $A_1u:=\partial^2_xu$, $A_2(u):=u(1-u)$. These operators are independent of time. Following the convention for linear operators we neglect the parenthesis in $A_1 u,$ at the same time we would like to emphasize that $A_2$ is nonlinear: $A_2(u)$.

\begin{center}
\begin{figure}[!hbt]
\includegraphics[width=8cm]{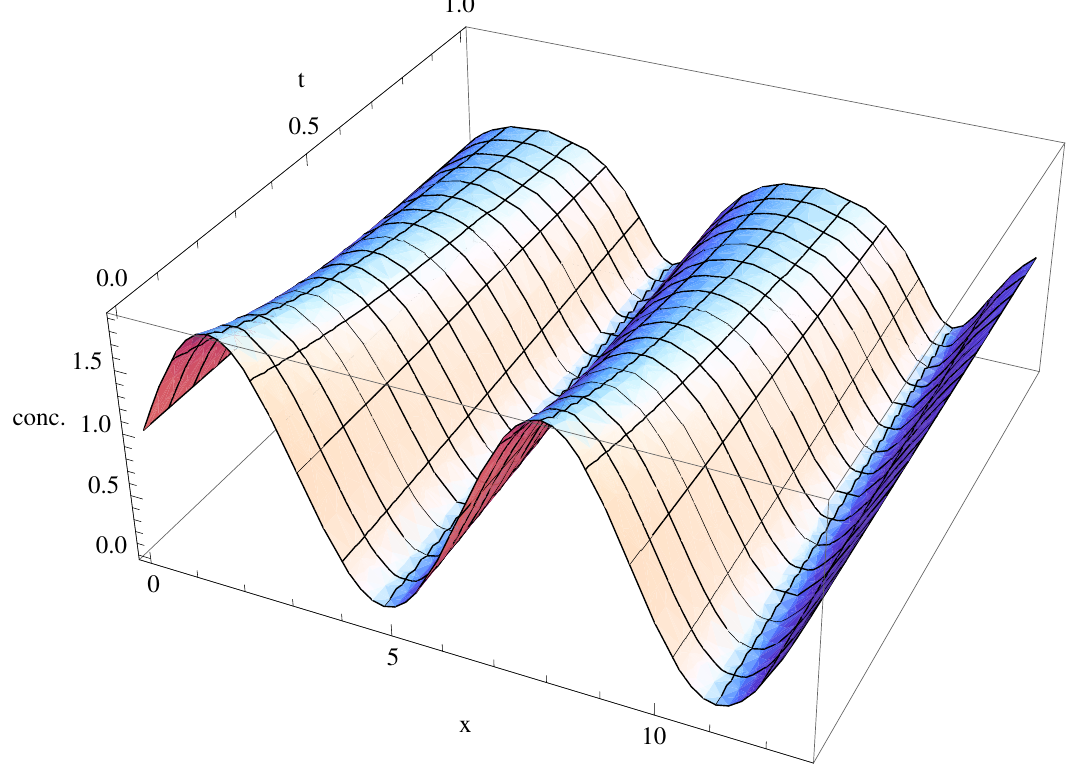}
\caption{The exact solution of (\ref{subprb}), $\displaystyle{\eta_2(x)=\frac{9}{10}\sin (x)+1,}$ $t\in [0,1]$ and $x\in [0,4\pi ].$ The same initial condition is used as in section 6.
}\label{regularfig1}
\end{figure}
\end{center}

The exact solution of problem (\ref{subprb}) is known, it is:
\begin{equation}\label{reac}
u_2(t,x)=\frac{\eta_2(x)e^t}{1-\eta_2(x)+\eta_2(x)e^t}.
\end{equation}
Since
\[
\lim_{t\to \infty} u_2(t,x)=\left\{
\begin{array}{lll}
1 & \text{when} & \eta_2(x) \neq 0\\
0 & \text{when} & \eta_2(x)= 0
\end{array}
\right.
\]
the solution has two stationary states, namely: $u_2(t,x)\equiv 0$,
$u_2(t,x)\equiv 1$. The $u_2(t,x)\equiv 1$ solution is asymptotically stable, whereas zero is an unstable equilibrium.

Knowing the exact solution of this subproblem as a function of the
initial condition means that we can symbolically solve this subproblem in each
time step during the splitting procedure. It might be worth using the
exact solution for comparisons in the study of the effect of
splitting. The exact solution of \eqref{subpra} is of no real use.

\section{Commutation of diffusion and reaction}
In this section we investigate the conditions under which the SEQ  has zero splitting error for reaction-diffusion systems in general. In reaction-diffusion equations there are two operators present on the right hand side, a linear and a nonlinear operator. Let us consider

\begin{equation}\label{rd}
U^\prime(t)=A(U(t))+R(U(t)) \quad U(0)=U_0,
\end{equation}
where $U:\mathbb{R}^+_0\to X$ with $X=\{\varphi=(\varphi_1,...,\varphi_M);\, \varphi_m: \mathbb{R}^3\to\mathbb{R}^+_0, \, m=1,...,M\}.$ We have $M$ different species undergoing the processes of diffusion and reactions. The operator $A$ represents the diffusion assuming that there is no cross diffusion present. In the above vectorial form $A$ is an operator matrix essentially with the spatial Laplacian in the diagonal and zeros everywhere else. $A$ is a linear operator. The operator $R$ usually acts as compositions with the multivariable polynomials $R_1, R_2...,R_M,$ it is a nonlinear operator. That is 
\[
R(\varphi)=(R_1(\varphi_1,...,\varphi_M),\dots,R_M(\varphi_1,...,\varphi_M))
\]
and
\[
A(\varphi)=A(\varphi_1,...,\varphi_M)=(D_1 \triangle \varphi_1,\dots, D_M \triangle \varphi_M).
\]
The global equation \eqref{rd} in the local form is:
\[
U^\prime(t)(x)=A(U(t))(x)+R(U(t))(x) \quad U(0)(x)=U_0(x),
\]
Therefore, another form of \eqref{rd} is
\begin{equation}
\left\{
\begin{array}{lcl}
\partial_t u_1(t,\textbf{x}) & = & D_1 \triangle u_1(t,\textbf{x}) +R_1(\ub(t,\textbf{x}))\\
& \vdots & \\
\partial_t u_M(t,\textbf{x}) & = & D_M \triangle u_M(t,\textbf{x}) +R_M(\ub(t,\textbf{x}))
\end{array}
\right.
\end{equation}
with $D_m>0,$ $m=1,\dots,M$, $\ub(\cdot)=(u_1(\cdot),...,u_M(\cdot))$ and $\textbf{x}\in \mathbb{R}^3.$ Since $A$ is linear its derivative is $A$ itself. The derivative of $R:=(R_1,...,R_M)$ which is the Jacobi matrix:
\[
R'=
\left(
\begin{array}{lcl}
\partial_1R_1 & \dots & \partial_M R_1\\
\vdots & & \vdots \\
\partial_1R_M & \dots & \partial_M R_M\\
\end{array}
\right)
\]
In the following we derive the condition of zero splitting error.
The sufficient condition of zero splitting error is \cite{verwer} that for every $\varphi \in X$ 
\begin{equation}\label{ver}
\left(A'\circ R-R'\circ A\right)(\varphi)=0
\end{equation}
Upon applying this we get
\[
0=D_m\triangle R_m(\varphi(t,\mathbf{x})) -\sum_{k=1}^M\partial_k R_m(\varphi(t,\mathbf{x})) D_m\triangle \varphi_k(t,\mathbf{x})
\]
For the sake of simplicity we will not carry the argument $(t,\mathbf{x}),$ it will no lead to misconceptions.
\[
0=\partial_x\sum_{k=1}^M\partial_kR_m(\varphi) \partial_x\varphi_k+\partial_y\sum_{k=1}^M\partial_k R_m(\varphi) \partial_y\varphi_k+\partial_z
\sum_{k=1}^M\partial_kR_m(\varphi) \partial_z\varphi_k-
\]

\[
-\sum_{k=1}^M\partial_k R_m(\varphi) \triangle \varphi_k \qquad \qquad m=1,...,M.
\]
Each of the first three terms is a sum of products thus, according to rule of differentiation of a product one gets:
\[
\sum_{k=1}^M\partial_x(\partial_k R_m(\varphi)) \partial_x \varphi_k+\sum_{k=1}^M\partial_y(\partial_k R_m(\varphi)) \partial_y\varphi_k+\sum_{k=1}^M\partial_z(\partial_k R_m(\varphi)) \partial_z\varphi_k
+
\]

\[
+\sum_{k=1}^M\partial_k R_m(\varphi) \partial_x^2\varphi_k+\sum_{k=1}^M\partial_k R_m(\varphi) \partial_y^2\varphi_k+\sum_{k=1}^M\partial_k R_m(\varphi) \partial_z^2\varphi_k
 -\sum_{k=1}^M\partial_k R_m(\varphi) \triangle \varphi_k
\]
Since the 4th, 5th, 6th and 7th terms eliminate each other, only the first three terms remain. Upon performing the differentiation the above expression becomes:
\[
\sum_{k=1}^M\sum_{j=1}^M\partial_j\partial_k R_m(\varphi) \partial_x\varphi_j \partial_x\varphi_k+\sum_{k=1}^M\sum_{j=1}^M\partial_j\partial_k R_m(\varphi) \partial_y\varphi_j \partial_y\varphi_k+
\]
\[
+\sum_{k=1}^M\sum_{j=1}^M\partial_j\partial_k R_m(\varphi) \partial_z\varphi_j \partial_z\varphi_k
\]
Using the notation:  $\partial_x\varphi:=(\partial_x\varphi_1,...,,\partial_x\varphi_M)$  we can reformulate this as:

\[
(R_m''(\varphi)\cdot\partial_x\varphi)\cdot\partial_x\varphi+(R_m''(\varphi)\cdot\partial_y\varphi)\cdot\partial_y\varphi
+(R_m''(\varphi)\cdot\partial_z\varphi)\cdot\partial_z\varphi=
\]

\[
=<R_m''(\varphi)\cdot\partial_x\varphi,\partial_x\varphi>+<R_m''(\varphi)\cdot\partial_y\varphi,\partial_y\varphi>
+<R_m''(\varphi)\cdot\partial_z\varphi,\partial_z\varphi>.
\]
The splitting error is zero if the above expression equals zero.
Using the notation:
\[
\partial\varphi:=(\partial_x\varphi,\partial_y\varphi,\partial_z\varphi)=\left(
\begin{array}{ccc}
\partial_x\varphi_1 & \partial_y\varphi_1 & \partial_z\varphi_1\\
\vdots&\vdots&\vdots \\
\partial_x\varphi_M & \partial_y\varphi_M & \partial_z\varphi_M
\end{array}
\right)
\]
we can write the formula above in the short form:
\[
(R_m''(\varphi)\cdot
\partial\varphi)^T
\cdot
\partial\varphi=0.
\]
Now we can formulate a statement.
\begin{thm} With the notations above the error of SEQ is zero if
for every function $\varphi$ and for all $m=1,2,...,M$  $(R_m''(\varphi)\cdot \partial \varphi)^T\cdot \partial \varphi=0$ holds.
\end{thm}

\begin{rmk}
We get the necessary condition of zero splitting error if we require condition (\ref{ver}) to hold only for the solution function of (\ref{rd}). Naturally without knowing the exact solution we can not check whether this holds or not. But the above formula provides us a sufficient condition for zero splitting error. If the equation holds for every possible function $\varphi$ then it will hold for the solution as well. This holds for every function iff all the entries of $R_m''$ is zero, which means that $R_m$ is a polynomial of at most first degree, $m=1,...,M:$ we only have first order reactions. Most of the practical problems have reaction terms of higher order, therefore there is almost always a splitting error. Our aim here is to examine the effect of splitting error in combined methods.
\end{rmk}

\begin{rmk}
Condition (\ref{ver}) is sufficient in the case of MS and SW splittings as well. Without going into the details condition \eqref{ver} ensures that $e^{\tau \mathbf{A}}e^{\tau \mathbf{R}}=e^{\tau (\mathbf{A}+\mathbf{R})},$ where $\mathbf{A}$ and $\mathbf{R}$ are the Lie-operators of $A$ and $R.$ In the case of MS we need $e^{\frac{\tau}{2} \mathbf{A}}e^{\tau \mathbf{R}}e^{\frac{\tau}{2} \mathbf{A}}=e^{\tau (\mathbf{A}+\mathbf{R})}$ which obviously holds if $e^{\tau \mathbf{A}}e^{\tau \mathbf{R}}=e^{\tau (\mathbf{A}+\mathbf{R})}$. Theorem $1$ remains valid in the case of MS and SW splittings.
\end{rmk}

\section{Order of combined methods}

When we solve partial differential equations we can use some kind of splitting but we can not avoid applying some numerical method as well. So in practice we use a combined method, a mixture of operator splitting and a numerical scheme and generate a solution for a nonlinear \PDE\ like (\ref{fisher}). We use the Taylor-formula to determine the order of the local error of this combined method. The Taylor-formula in normed vector spaces can be found in e. g. \cite{komornik}. Here we recall the Taylor-formula for normed vector spaces:
\begin{thm} If $f:X\to Y$ is $n$ times differentiable in $a\in X$ and $h\rightarrow0,$ then
\[f(a+h)=\sum_{k=0}^n \frac{f^{(k)}(a)}{k!}h^k+\epsilon(h) \|h\|^n\]
where $h^k:=(h,...,h)\in X^k$ and $\displaystyle{\lim_{h\to 0} \varepsilon (h)=0.}$
\end{thm}
Suppose that $U$ is the solution of the equation:
\begin{equation}\label{uder}
U^{\prime}(t) = AU(t)+R(U(t))
\end{equation}
From now on in this section we restrict our investigations to problems where $A:X\to X$ is a bounded linear operator defined on the whole set $X,$ therefore the following derivations are not directly applicable to the case of diffusion. $A$ is differentiable and its derivative is $A$ itself for every $x\in X.$ The operator $R:X\to X$ will act as a composition with a differentiable nonlinear function thus $R$ is a differentiable mapping as well. Based on \eqref{uder} and the chain-rule $U^\prime$ is a differentiable function and $U^\prime(t)\in X.$

The Taylor-expansion of $U$ in time $t_0$ is 
\[
U_\tau:=U(t_0+\tau)=U(t_0)+U^{\prime}(t_0) \tau+\frac{1}{2}U''(t_0)\tau^2+\varepsilon(\tau)  \| \tau \|^2.
\]
The norm we will neglect from now on since $U:\mathbb{R}\to X,$ $\tau$ denotes a positive real number.
$U'(t_0)$ is given by \eqref{uder}, we get $U''(t_0)$ by differentiation of \eqref{uder}: 
\[
U'' (t_0)=A'(U(t_0))\circ U^\prime (t_0)+R'(U(t_0))U^\prime (t_0)=
A(U^\prime (t_0))+R'(U(t_0))U^\prime (t_0)=
\] 
\[
=A(A(U(t_0))+R(U(t_0)))+R'(U(t_0))A(U(t_0))+R'(U(t_0))R(U(t_0)).
\]
$U(t_0)$ is denoted by $U_0.$ Using this notation if the value of $U(t_0)=U_0$ is known then we can approximate $U_\tau:$
\begin{equation}\label{exutau}
U_\tau=U_0+\big(A(U_0)+R(U_0)\big) \tau+
\end{equation}
\[
+\frac{1}{2}\Big( A(A(U_0))+A(R(U_0))+R'(U_0)A(U_0)+R'(U_0)R(U_0)\Big)\tau^2+\varepsilon(\tau) \tau^2.
\]
Beyond the conditions we already mentioned this form of $U_\tau$ exists if $U_0\in D(A^2)$ and $R(U_0)\in D(A)$ which naturally hold if $A:X\to X$ and $R:X\to X.$ For the theorems in this section we will need $R$ to be three times continuously differentiable.
\subsection{Methods of first order}
\begin{thm}
The sequential splitting combined with the first order Euler forward scheme provides a first order method.
\end{thm}
The proof will be given in two steps.

\subsubsection{Linear-Nonlinear}

If we use SEQ starting with the nonlinear problem corresponding to $R$ combined with Euler forward method for both subproblems we get:

\begin{equation*}
\left\{ {\setlength\arraycolsep{.13889em}
\begin{array}{rll}
V(\tau)& = &U_0+ \tau R(U_0)\\
\bar{U}_\tau & = & V(\tau)+\tau A(V(\tau)) \\
\end{array}}
\right.
\end{equation*}

\Pf{
The approximation of the solution in time $\tau$ is:
\[
\tilde{U}_\tau=V(\tau)+\tau A(V(\tau))=U_0+ \tau R(U_0)+\tau A(U_0+ \tau R(U_0)).
\]
Since $A$ is linear we have:
\begin{equation}\label{apprutau}
\tilde{U}_\tau=U_0+ \tau R(U_0)+\tau A(U_0)+ \tau^2 A(R(U_0)).
\end{equation}
The local error generated in this step of length $\tau$ based on \eqref{exutau} and \eqref{apprutau} is:
\[
U_\tau-\tilde{U}_\tau=\Big(A(A(U_0))-A(R(U_0))+R'(U_0)A(U_0)+R'(U_0)R(U_0)\Big)\frac{\tau ^2 }{2}+\varepsilon(\tau) \tau^2.
\]
}

\subsubsection{Nonlinear-Linear}

If we use SEQ starting with the linear problem corresponding to $A$ combined with Euler forward method for both subproblems we get:

\begin{equation*}
\left\{ {\setlength\arraycolsep{.13889em}
\begin{array}{rll}
V(\tau)& = &U_0+ \tau A(U_0)\\
\bar{U}_\tau & = & V(\tau)+\tau R(V(\tau)). \\
\end{array}}
\right.
\end{equation*}

\Pf{
The approximation of the solution in time $\tau$ is:
\[
\tilde{U}_\tau=V(\tau)+\tau R(V(\tau))=U_0+ \tau A(U_0)+\tau R(U_0+ \tau A(U_0)).
\]
Let us define the function $F:\mathbb{R}\to X$ in the following way: $F(\delta):=R(U_0+\delta \tau A(U_0))$. Then
\[
F(0)=R(U_0) \text{ and } F(1)=R(U_0+\tau A(U_0)),
\]
and $F$ is differentiable since it equals to $R\circ f$ with $f(\delta):=U_0+\delta \tau A(U_0)$ (where is differentiable), and according to the chain-rule 
\[
F'(\delta)=R'(U_0+\delta \tau A(U_0)) \tau A(U_0) \quad F''(\delta)=R''(U_0+\delta \tau A(U_0)) (\tau A(U_0))^2
\]
and
\[
F^{(n)}(\delta)=R^{(n)}(U_0+\delta \tau A(U_0))(\tau A(U_0))^n
\]
For the Taylor-expansion of $F$ we need a similar but more specific relation.
\begin{prp} If $F:\mathbb{R}\to X$ is $n$ times differentiable in every point of $[a, b]$ then there exists $c\in (a,b)$ such that 
\[
\|F(b)-\sum_{k=0}^{n-1} \frac{F^{(k)}(a)}{k!}\| \leqslant \frac{\|F^{(n)}(c)\|}{n!} (b-a)^n.
\]
\end{prp}
In other words 
\[
F(b)=\sum_{k=0}^{n-1} \frac{F^{(k)}(a)}{k!}+ \varepsilon_1(c) \frac{\|F^{(n)}(c)\|}{n!} (b-a)^n
\]
with $\|\varepsilon_1(c)\| \leqslant 1.$
The expansion of $F$ around $0$ gives:
\[
F(1)=F(0)+F'(0)+\varepsilon_1(c)\frac{1}{2}\|F''(c)\| \qquad c\in[0,1]
\]
implying
\[
R(U_0+\tau A(U_0))=R(U_0)+R'(U_0) \tau A(U_0)+\varepsilon_1(c)\frac{1}{2}\|R''(U_0+c\tau A(U_0)) (\tau A(U_0))^2\|=
\]
\[
=R(U_0)+R'(U_0) \tau A(U_0)+\varepsilon(\tau)\tau,
\]
where
\[
\varepsilon(\tau)=\varepsilon_1(c)\frac{1}{2}\|R''(U_0+c\tau A(U_0)) A(U_0)^2\tau \|.
\]
This tends to zero if $\tau$ tends to zero. This is the only relevant property of $\varepsilon$ here. Although $c$ can change as $\tau$ changes but since $\|\varepsilon_1(c)\| \leqslant 1$ we can ignore $\varepsilon$'s dependence on $c$ through $\varepsilon_1.$
\[
\tilde{U}_\tau=U_0+ \tau A(U_0)+\tau \left(R(U_0)+R'(U_0) \tau A(U_0)+\varepsilon(\tau)  \| \tau \|\right)=
\]
\[
=U_0+ \tau A(U_0)+\tau R(U_0)+R'(U_0) \tau^2 A(U_0)+\varepsilon(\tau) \tau^2,
\]
\[
U_\tau=U_0+ \tau (A(U_0)+R(U_0))+
\]
\[
+\Big(A^2(U_0)+A(R(U_0))+R'(U_0)A(U_0)+R'(U_0)R(U_0)\Big)\frac{\tau ^2}{2}+\varepsilon(\tau) \tau^2.
\]
Comparing this with the approximation we get: 
\[
U_\tau-\tilde{U}_\tau=\Big(A^2(U_0)+A(R(U_0))-R'(U_0)A(U_0)+R'(U_0)R(U_0)\Big)\frac{\tau ^2 }{2}+\varepsilon(\tau) \tau^2.
\]
}

\subsubsection{Weighted Splitting}

\begin{thm} SW combined with the first order Euler forward method provides a method of first order.
\end{thm}
\Pf{
Here we simply use the above results with some $\omega\in[0,1]$ parameter: 
\[
\tilde{U}_\tau=\omega \big(U_0+ \tau \left(R(U_0)+A(U_0)\right)+ \tau^2 A(R(U_0))\big)+\]
\[
+(1-\omega)\big(U_0+ \tau \left(A(U_0)+R(U_0)\right)+R'(U_0) \tau^2 A(U_0)+R''(U_0) \tau^3 (A(U_0))^2\frac{1}{2}+\varepsilon(\tau) \tau^2,
\]
\[
\tilde{U}_\tau=U_0+ \tau \left(A(U_0)+R(U_0)\right)+\big(\omega A(R(U_0))+(1-\omega)(R'(U_0) A(U_0))\big) \tau^2+\varepsilon(\tau) \tau^2.
\]
For the local error we have:
\[
U_\tau-\tilde{U}_\tau=
\]
\[
=\Big(A^2(U_0)+(1-2\omega)A(U_0)R(U_0)+(2\omega-1)R'(U_0)A(U_0)+R'(U_0)R(U_0)\Big)\frac{\tau ^2 }{2}+
\]
\[
+\varepsilon(\tau) \tau^2,
\]
for $\displaystyle{\omega=\frac{1}{2}}$ we have 
\[
U_\tau-\tilde{U}_\tau=\Big(A^2(U_0)+R'(U_0)R(U_0)\Big)\frac{\tau ^2 }{2}+\varepsilon(\tau) \tau^2.
\]
}
Conclusion: Although the SW is of second order its combination with the first order Euler method provides only first order accuracy.

\subsection{Methods of higher order}

The above derivation can be used to determine the order of combined methods with higher order numerical schemes.
The method can be extended for schemes of arbitrary order although the calculations become very complicated as the order increases. As an example let us consider the improved Euler scheme which is of second order and combine it with SEQ:

\begin{thm} The second order improved Euler scheme combined with SEQ provides a first order method.
\end{thm}
Again, the proof will be given in two steps.

\subsubsection{Nonlinear-Linear}

\begin{equation*}
\left\{ {\setlength\arraycolsep{.13889em}
\begin{array}{rll}
V(\tau)& = &U_0+ \tau A(U_0+\displaystyle{\frac{\tau}{2}A(U_0))}\\
\bar{U}_\tau & = & V(\tau)+\tau R(V(\tau)+\displaystyle{\frac{\tau}{2}R(V(\tau)))} \\
\end{array}}
\right.
\end{equation*}

\Pf{
The approximation of the solution is in time $\tau$ is:
\[
\tilde{U}_\tau=V(\tau)+\tau R\Big(V(\tau)+\frac{\tau}{2}R(V(\tau))\Big)=U_0+ \tau A\big(U_0+\frac{\tau}{2}A(U_0)\big)+
\]
\[
+\tau R\Big(U_0+ \tau A\big(U_0+\frac{\tau}{2}A(U_0)\big)+ \frac{\tau}{2}R\Big(U_0+ \tau A\big(U_0+\frac{\tau}{2}A(U_0)\big)\Big)\Big)=
\]
\[
=U_0+ \tau A(U_0)+\frac{\tau^2}{2}A^2(U_0)+
\]
\[
+\tau R\Big(U_0+ \underline{\tau \Big(A(U_0)+\frac{\tau}{2}A^2(U_0)+ \frac{1}{2}R\big(U_0+\tau A(U_0)+\frac{\tau^2}{2}A^2(U_0)\big)\Big)}\Big)=
\]
The underlined part is the coefficient of $\delta$ in the argument of $F$. The first order Taylor--expansion gives:
\[
=U_0+ \tau A(U_0)+\frac{\tau^2}{2}A^2(U_0)+\tau R(U_0)+
\]
\[+\tau R'(U_0) \underline{\tau \Big(A(U_0)+\frac{\tau}{2}A^2(U_0)+ \frac{1}{2}R\big(U_0+\tau A(U_0)+\frac{\tau^2}{2}A^2(U_0)\big)\Big)}+\varepsilon(\tau) \tau^2=
\]
\[
=U_0+ \tau \big(A(U_0)+R(U_0)\big)+
\]
\[
+\frac{\tau^2}{2}\Big(A^2(U_0)+2R'(U_0) A(U_0)+R'(U_0)R\big(U_0+\underline{\tau A(U_0)+\frac{\tau^2}{2}A^2(U_0)}\big)\Big)+
\]
\[
+\frac{\tau^3}{2}R'(U_0)A^2(U_0)+\varepsilon(\tau) \tau^2=
\]
Taking the Taylor--expansion again, the coefficient of $\delta$ is the underlined part:
\[
=U_0+ \tau \big(A(U_0)+R(U_0)\big)+
\]
\[
+\frac{\tau^2}{2}\Big(A^2(U_0)+2R'(U_0) A(U_0)+ R'(U_0)\big(R(U_0)+\tau R'(U_0)A(U_0)+\frac{\tau^2}{2}R'(U_0)A^2(U_0))\big)\Big)+
\]
\[
+\varepsilon(\tau) \tau^2=
\]
\[
=U_0+ \tau \big(A(U_0)+R(U_0)\big)+\frac{\tau^2}{2}\Big(A^2(U_0)+2R'(U_0) A(U_0)+ R'(U_0)R(U_0)\Big)+\varepsilon(\tau) \tau^2
\]
\[
U_\tau-\tilde{U}_\tau=\Big(A(R(U_0))-R'(U_0)A(U_0)\Big)\frac{\tau ^2 }{2}+\varepsilon(\tau) \tau^2.
\]
}

\subsubsection{Linear-Nonlinear} The proof of the linear-nonlinear case is more straightforward:

\Pf{
\begin{equation*}
\left\{ {\setlength\arraycolsep{.13889em}
\begin{array}{rll}
V(\tau)& = &U_0+ \tau R(U_0+\displaystyle{\frac{\tau}{2}R(U_0))}\\
\tilde{U}_\tau & = & V(\tau)+\tau A(V(\tau)+\displaystyle{\frac{\tau}{2}A(V(\tau)))} \\
\end{array}}
\right.
\end{equation*}
\[
\tilde{U}_\tau=V(\tau)+\tau A\Big(V(\tau)+\frac{\tau}{2}A(V(\tau))\Big)=
\]
\[
=U_0+ \tau R\big(U_0+\frac{\tau}{2}R(U_0)\big)+
\]
\[
+
\tau A\Big(U_0+ \tau R\big(U_0+\frac{\tau}{2}R(U_0)\big)+ \frac{\tau}{2}A\Big(U_0+ \tau R\big(U_0+\frac{\tau}{2}R(U_0)\big)\Big)\Big)=
\]
\[
=U_0+ \tau \big(R(U_0)+\frac{\tau}{2}R(U_0)R'(U_0)\big)+\tau A(U_0)+
\]
\[
+\tau^2\Big(A\big( R(U_0)+\frac{\tau}{2}R(U_0)R'(U_0)\big)+ \frac{1}{2}A^2(U_0)\Big)+\varepsilon(\tau) \tau^2=
\]
\[
=U_0+ \tau \big(R(U_0)+A(U_0)\big)+\frac{\tau^2}{2}\Big(R(U_0)R'(U_0)+2A(R(U_0))+A^2(U_0)\Big)+\varepsilon(\tau) \tau^2
\]
\[
U_{\tau}-\tilde{U}_\tau=\Big(R'(U_0)A(U_0)-A(R(U_0))\Big)\frac{\tau^2}{2}+\varepsilon(\tau) \tau^2.
\]
}

As we can see this combined method is of first order. Although the applied numerical scheme ensures second order accuracy the use of sequential splitting results in order reduction. We followed the same ideas and calculated the orders for combinations of the introduced splitting methods and four different numerical schemes. The table below contains our results on orders of different splittings coupled with different numerical methods. Symbolic calculations on for example MS splitting coupled with 4th order Runge--Kutta method becomes complicated. An algorithm was written in \textit{Mathematica} for these symbolic calculations.

\vspace{8mm}
\begin{table}[h]
\begin{tabular}{|c||c|c|c|c|}
\hline
   & exp. Euler (1) & impr. Euler (2) & Heun (3)& Runge--Kutta (4)\\
\hline \hline
 SEQ  $(1)$ & $1$ & $1$ & $1$ & $1$\\
\hline
 SW  $(2)$ & $1$ & $2$  & $2$ & $2$\\
\hline
 MS  $(2)$ & $1$ & $2$ & $2$ & $2$\\
 \hline
\end{tabular}
\caption{Local orders of  combined methods for (\ref{uder})}
\end{table}

\vspace{8mm}

The order of the methods are in the parenthesis. A study of the order of combined methods for bounded linear problems can be found in \cite{csomosfarago}. Their results say that the order of the combined method is $s:=\min\{p,r\},$ if $p$ is the order of splitting and $r$ denotes the order of the numerical method. The numbers of the above table are in accordance with their results.

\section{Numerical experiments}

Here we introduce our numerical results on the Fisher equation. We solved both subproblems (\ref{subpra}), (\ref{subprb}) using numerical methods of four different orders; the explicit Euler, the improved Euler method which is of second order, the third order Heun and the fourth order Runge--Kutta method. We investigate the SEQ, SW and MS splitting methods which are of first and second orders. We calculated the errors and orders of these combined methods numerically. Our test problem is the following initial--boundary value problem:

\begin{equation}\label{num}
\left\{ {\setlength\arraycolsep{.13889em}
\begin{array}{lll}
\partial_tu (t,x)& = & \partial^2_xu(t,x)+u(t,x)(1-u(t,x)) \\
u(0,x) & = & 1+0.9\sin(x) \\
u(t,0) & = & 1 \\
u(t,4\pi) & = & 1
\end{array}}
\right.,
\end{equation}
where $x\in [0, 4\pi]$ and $t\in [0,1]$.

\begin{center}
\begin{figure}[!hbt]
\includegraphics[width=8cm]{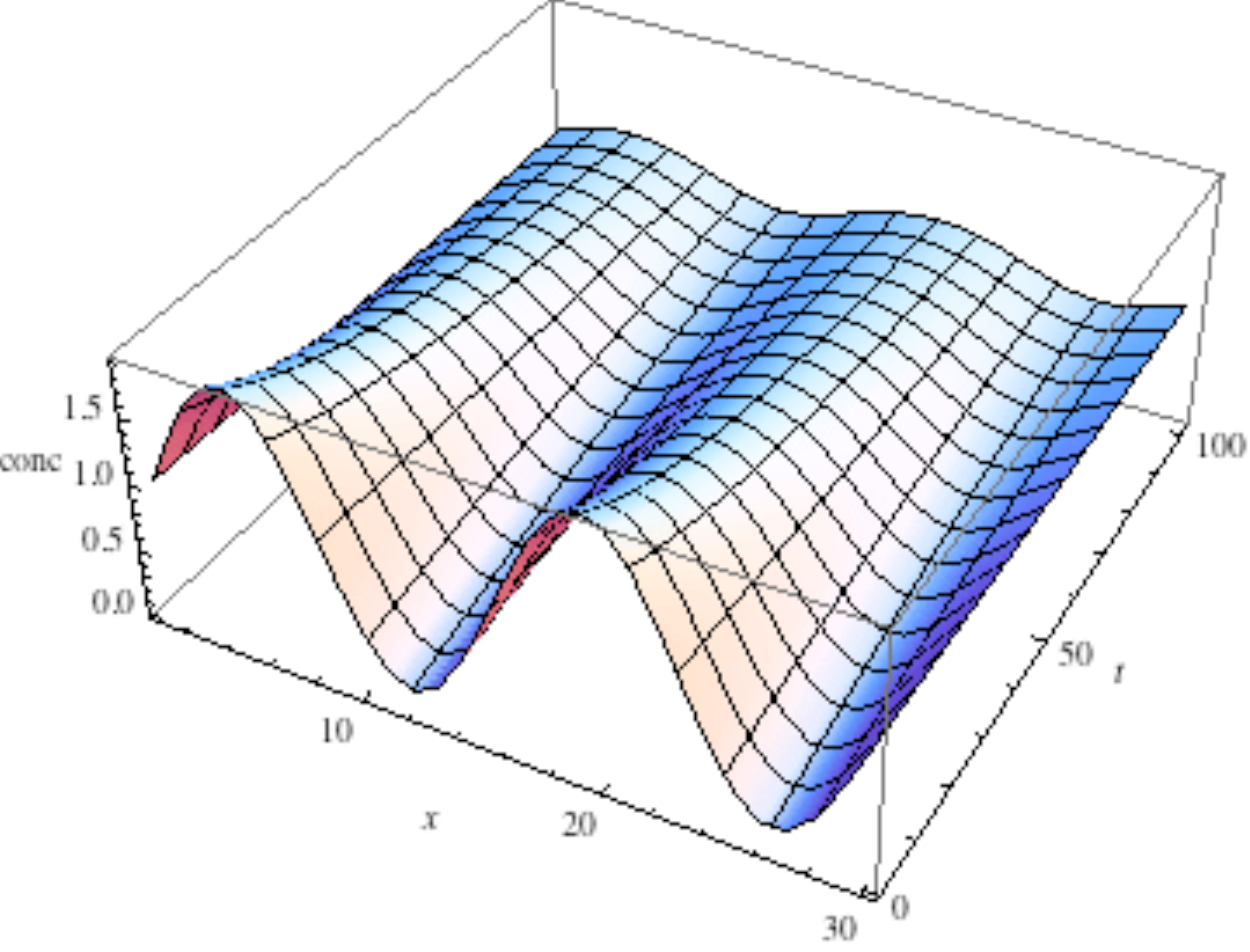}
\caption{Reference solution generated by fourth order Runge--Kutta scheme, $\tau=0.01.$}
\label{regularfig1}
\end{figure}
\end{center}

We performed a spatial semidiscretization with length parameter $\Delta x=\frac{4\pi}{30}$ that is we divided $[0,4\pi]$ into $N=30$ parts of equal length. Our tests showed that finer divisions provides no significantly more accurate solutions that is the obtained error is of the same magnitude as with  $N=30.$ We approximated the spatial derivative with the well known second order scheme:
\[
\partial^2_x u(t,x_i)\approx \frac{u(t,x_{i+1})-2u(t,x_i)+u(t,x_{i-1}))}{\Delta x^2}.
\]
 After temporal discretization with $\tau=0.01$ we solved the full problem \eqref{num} with the fourth order Runge--Kutta method. Taking a smaller time step resulted in solution that differs only in a magnitude of $10^{-6}$. This provides the reference solution for our study. In the experiments we used spatial division $N=30$ in every case, in fact we investigated the convergence of the semidiscrete submodels to the semidiscrete model: the reference solution. On the connection between the convergence to a semidiscrete model and convergence to a continuous model see \cite{thomee}.

\subsection{Determination of the local order}

The local error $E(\tau)$ of a method is of order $s$ if
\[
E(\tau)=O(\tau^{s+1})
\]
or by using a different formulation the order is $s$ if
\[
\lim_{\tau\rightarrow0}\frac{E(\tau)}{\tau^{s+1}}
\]
is finite where $s\in\mathbb{N}$ is the smallest number with this
property. Thus in practice we can estimate the order of the
combined method by calculating this limit with different fixed values of
$s$ until we find the appropriate one. For a fixed numerical
step size $h$ the limit means naturally that $\tau\rightarrow h,$
therefore we choose $h$ to be as small as possible close to the
smallest number that our computer can represent.

Another way to calculate the order of the method is the following.
Considering the two formula above we can conclude that
\[
\frac{E(\tau)}{\tau^{s+1}}\approx c
\]
for small $\tau$-s, where $c$ is a constant which does not depend
on $\tau.$ We can take the logarithm of both sides:
\[
\log E(\tau)\approx(s+1)\log \tau +\log c.
\]
This defines a straight line whose steepness gives the order of the method. The table below contains results for the local order $s$. The order of splitting is $p$ and $r$ is the order of the numerical scheme. The time step of the numerical method $h$ was the tenth of splitting time step $\tau$.

\vspace{8mm}
\begin{table}[h]
\begin{tabular}{|c||c|c|c|}
\hline
 $s$ $(h=0.1\tau )$  & exp. Euler $(r=1)$ & Heun $(r=3)$& Runge--Kutta $(r=4)$\\
\hline \hline
 SEQ  $(p=1)$ & $0.98$ & $0.98$ & $0.98$\\
\hline
 SW  $(p=2)$ & $0.83$  & $1.99$ & $1.99$\\
\hline
 MS  $(p=2)$ & $0.93$ & $1.96$ & $1.96$\\
 \hline
\end{tabular}
\caption{Orders of  combined methods for (\ref{num})}
\end{table}
It is interesting that with the explicit Euler method increasing the order of splitting does not improve the results.

\subsection{Estimation of the global order}

We used the time steps $\tau_i$ in such a way that the evaluation of the error is simple, since the corresponding division of the time interval is a subset of the one of $\tau=0.01$.

\vspace{5mm}

\begin{tabular}{|c|c|c|c|c|c|c|}
\hline 
$\tau_1$& $\tau_2$ &$\tau_3$ & $\tau_4$ & $\tau_5$ & $\tau_6$ & $\tau_7$\\
 \hline
 0.2 & $0.1$ & $0.0625$ & $0.05$ & $0.04$ & $0.025$ &$0.02$\\
\hline
\end{tabular}

\vspace{10mm}

As for the error we know that:
\[
E(\tau_1)\approx c\cdot \tau_1^{\rho}
\]
for small $\tau$-s, where $c$ is a constant which does not depend
on $\tau.$ So we can write:
\[
E(\tau_2)\approx c\cdot \tau_2^{\rho}
\]

\[
\frac{E(\tau_1)}{E(\tau_2)}\approx \left(\frac{\tau_1}{\tau_2}\right)^{\rho}
\]
We can take the logarithm of both sides:
\[
\log \frac{E(\tau_1)}{E(\tau_2)}\approx \rho \log \frac{\tau_1}{\tau_2}.
\]
For each $i$:
\[
\log \frac{E(\tau_i)}{E(\tau_{i+1})}\approx \rho \log \frac{\tau_i}{\tau_{i+1}}.
\]
So 
\[
\frac{\log \frac{E(\tau_i)}{E(\tau_{i+1})}}{ \log \frac{\tau_i}{\tau_{i+1}}}\approx \rho.
\]
Evaluation of the left side shall give us the same value for every $i=1,2,3,4,5,6$. The following table contains the results of this calculation for different splittings and numerical methods.

\vspace{5mm}
\begin{table}[h]
\begin{tabular}{|c||c|c|c|c|}
\hline
$\rho$ $(h=\tau)$ &$r=1$ & $r=2$ & $r=3$ & $r=4$\\
\hline \hline
 SEQ  $(p=1)$ & $1.04$ & $0.99$ & $1.08$ & $1.08$\\
\hline
 SW  $(p=2)$ & $1.02$ & $2.07$  & $2.01$ & $1.98$\\
\hline
 MS  $(p=2)$ & $1.02$ & $2.07$ & $1.95$ & $1.998$\\
 \hline
\end{tabular}
\caption{Orders of  combined methods for (\ref{num})}
\end{table}
Here we had the same experience as with the calculations shown in Table 2. Since the solution of (\ref{subprb}) is given in (\ref{reac}) as a function of the initial condition we can use it in calculations instead of the numerical solution. The table below contains results generated by using (\ref{reac}) in each time step.

\begin{table}[h]
\begin{tabular}{|c||c|c|c|c|}
\hline
$\rho$ $(h=\tau)$ &$r=1$&$r=2$&$r=3$&$r=4$\\
\hline \hline
 SEQ  $(p=1)$&$1.03$&$1.02$&$1.01$&$1.01$\\
\hline
 SW  $(p=2)$ & $1.01$ & $2.06$ & $1.95$ & $1.98$\\
\hline
 MS $(p=2)$ & $1.03.$ & $2.00$ & $1.99$ & $1.99$\\
 \hline
\end{tabular}
\caption{Orders of  combined methods for (\ref{num}), using (\ref{reac})}
\end{table}

The first column fits into the general scheme.

\subsection{Splitting into three operators}

Above we saw that the splitting error is zero if the reaction term is of first order. Another, natural decomposition of the right hand side of the Fisher equation is when we separate the reaction part into two terms. The split problem is the following:
\torol{
\begin{tabular}{lcr}
\begin{equation}
\left\{ {\setlength\arraycolsep{.13889em}
\begin{array}{lll}
\partial_tu_1 (t,x)& = & \partial^2_xu_1 (t,x)\\
u_1(0,x) & = & \eta_1(x) \\
\end{array}}
\right.
\end{equation}
&
\begin{equation}
\left\{ {\setlength\arraycolsep{.13889em}
\begin{array}{lll}
\partial_tu_2 (t,x)& = & u_2(t,x)\\
u_2(0,x) & = & \eta_2(x) \\
\end{array}}
\right.
\end{equation}
&
\begin{equation}
\left\{ {\setlength\arraycolsep{.13889em}
\begin{array}{lll}
\partial_tu_3 (t,x)& = & -u^2_3(t,x) \\
u_3(0,x) & = & \eta_3(x) \\
\end{array}}
\right.
\end{equation}
\end{tabular}
}
\begin{equation}
\begin{array}{lcr}
\left\{ {\setlength\arraycolsep{.13889em}
\begin{array}{lll}
\partial_tu_1 (t,x)& = & \partial^2_xu_1 (t,x)\\
u_1(0,x) & = & \eta_1(x) \\
\end{array}}
\right.
&
\left\{ {\setlength\arraycolsep{.13889em}
\begin{array}{lll}
\partial_tu_2 (t,x)& = & u_2(t,x)\\
u_2(0,x) & = & \eta_2(x) \\
\end{array}}
\right.
&
\left\{ {\setlength\arraycolsep{.13889em}
\begin{array}{lll}
\partial_tu_3 (t,x)& = & -u^2_3(t,x) \\
u_3(0,x) & = & \eta_3(x) \\
\end{array}}
\right.
\end{array}
\end{equation}
Then $A_1u=\partial^2_x u,$ $A_2u=u,$ $A_3(u)=-u^2.$ Since $A_2$ is a first order polynomial the splitting error of the SEQ $S_1S_2$ is zero according to theorem 1. It is easy to prove, that it is not zero for $S_3$ coupled with any of the other two operators.

\begin{table}[h]
\begin{tabular}{|c||c|c|c|c|}
\hline
s $(h=\tau)$ &$r=1$&$r=2$&$r=3$&$r=4$\\
\hline \hline
 SEQ  $(p=1)$ & $0.96$ & $1.96$ & $2.97$ & $3.96$\\
\hline
 SW $(p=2)$ & $0.96$ & $1.96$ & $2.97$ & $3.96$\\
\hline
 MS $(p=2)$ & $0.96$ & $1.96$ & $2.96$ & $3.96$\\
 \hline
\end{tabular}
\caption{Orders for splitting of $\partial_t u=\partial^2_xu+u$}
\end{table}

In the solution of subproblems associated with operator $A_2$ and $A_3$ we used the exact solution. The first subproblem was solved numerically. Considering the MS type splittings $S_3S_2S_1S_2S_3$ and $S_2S_3S_1S_3S_2$ it is reasonable to expect more accurate solutions with in the case $S_3S_2S_1S_2S_3$ since the neighbors $S_1$ and $S_2$ generate no splitting error. Whereas $S_2S_3S_1S_3S_2$ generate splitting errors between every neighboring operators. Figure 5 and figure 6 proves that although both MS provides the expected second order accuracy (combined with a third order or fourth order numerical schemes) the hypothetic relation in accuracy turns out to be the opposite. 

\begin{center}
\begin{figure}[!hbt]
\includegraphics[width=8cm]{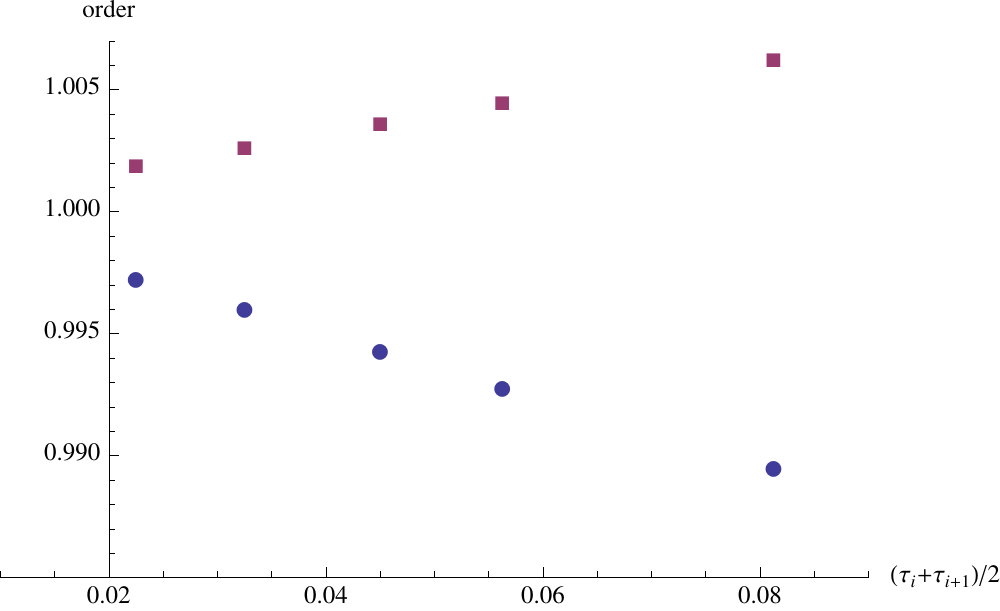}
\caption{Order approximation of the $S_1S_2S_3$ (squares) versus $S_1S_3S_2$ (discs) SEQ type splittings. Generated by the fourth order Runge-Kutta scheme with time steps $\tau=0.1,\, 0.0625,\, 0.05,\, 0.04,\, 0.025,\, 0.02$.
}\label{regularfig1}
\end{figure}
\end{center}

\begin{center}
\begin{figure}[!hbt]
\includegraphics[width=8cm]{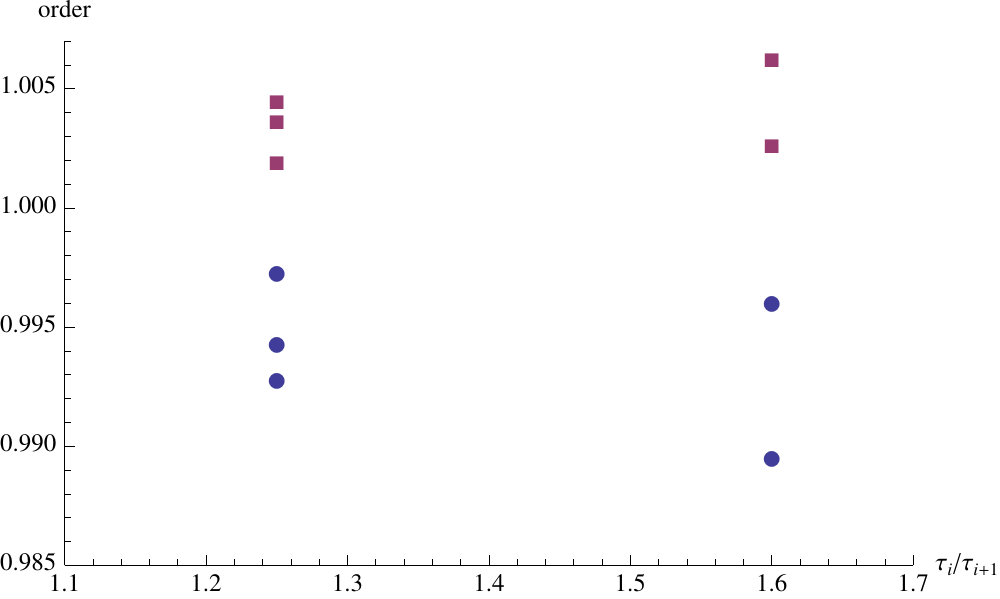}
\caption{Order approximation of SEQ type splittings as the function of $\tau_i/\tau_{i+1}.$
}\label{regularfig1}
\end{figure}
\end{center}

\begin{center}
\begin{figure}[!hbt]
\includegraphics[width=8cm]{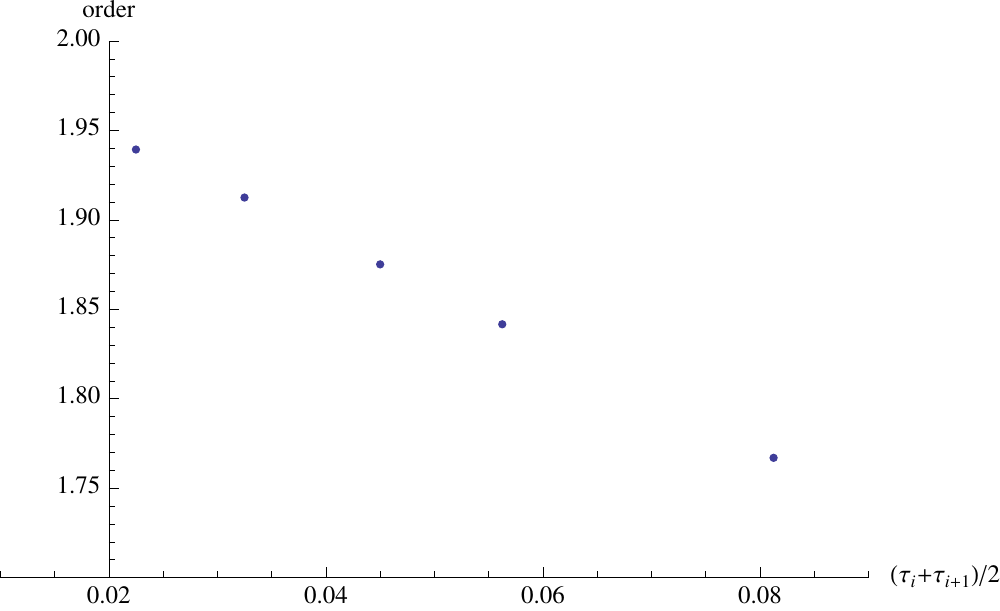}
\caption{Order approximation of the SW splitting which is the arithmetic mean of the six possible sequential splittings. Generated by the fourth order Runge-Kutta scheme.
}\label{regularfig1}
\end{figure}
\end{center}

\begin{center}
\begin{figure}[!hbt]
\includegraphics[width=8cm]{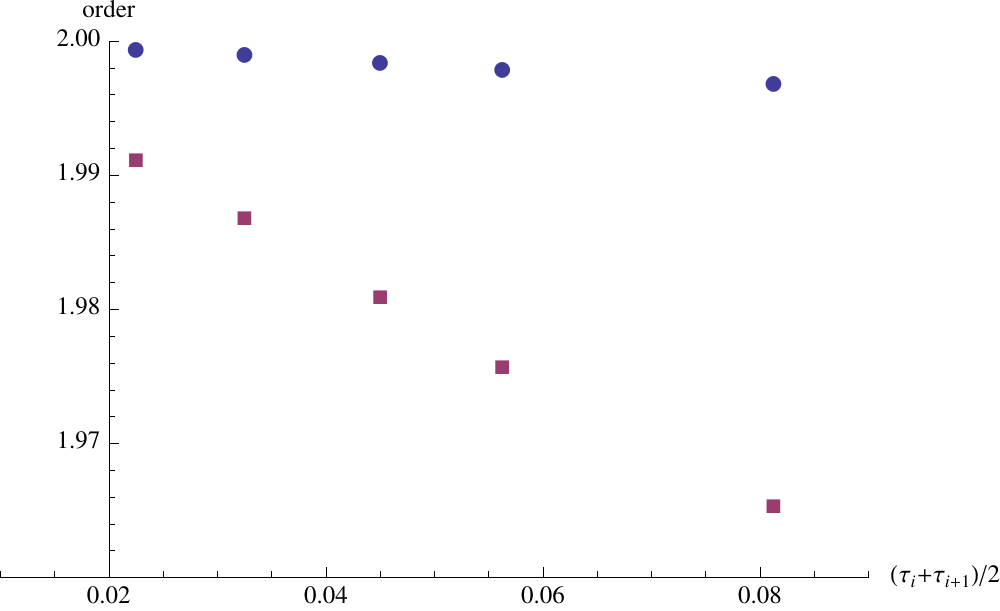}
\caption{Order approximation of the $S_3S_2S_1S_2S_3$ (discs) versus $S_2S_3S_1S_3S_2$ (squares) MS type splittings. Generated by the third order Heun scheme.
}\label{regularfig1}
\end{figure}
\end{center}

\begin{center}
\begin{figure}[!hbt]
\includegraphics[width=8cm]{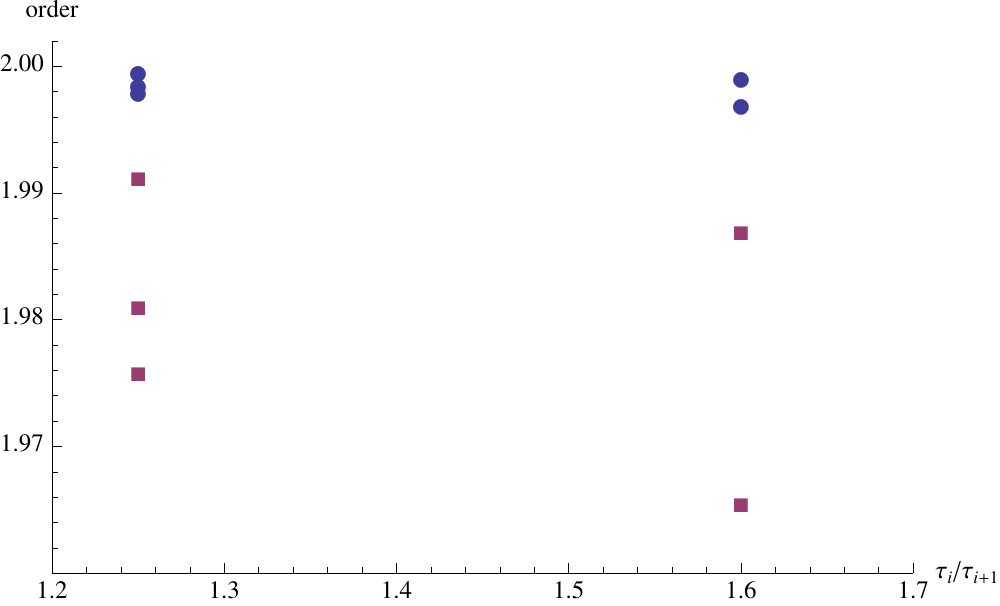}
\caption{Order approximation of MS type splittings as the function of $\tau_i/\tau_{i+1}.$
}\label{regularfig1}
\end{figure}
\end{center}

\begin{center}
\begin{figure}[!hbt]
\includegraphics[width=8cm]{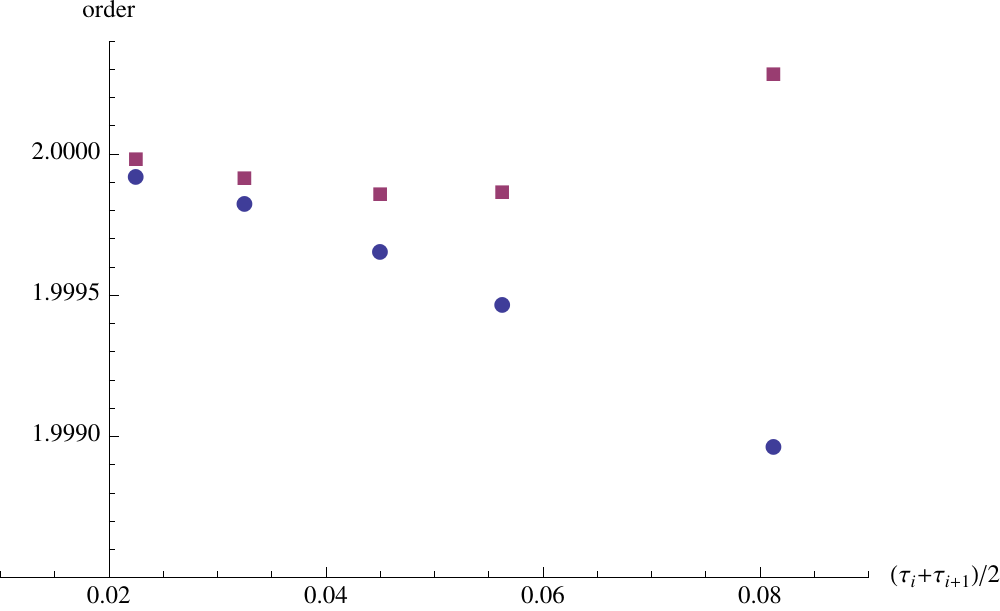}
\caption{Order approximation of the $S_3S_2S_1S_2S_3$ (discs) versus $S_2S_3S_1S_3S_2$ (squares) MS type splittings. Generated by the fourth order Runge-Kutta scheme.
}\label{regularfig1}
\end{figure}
\end{center}

\begin{center}
\begin{figure}[!hbt]
\includegraphics[width=8cm]{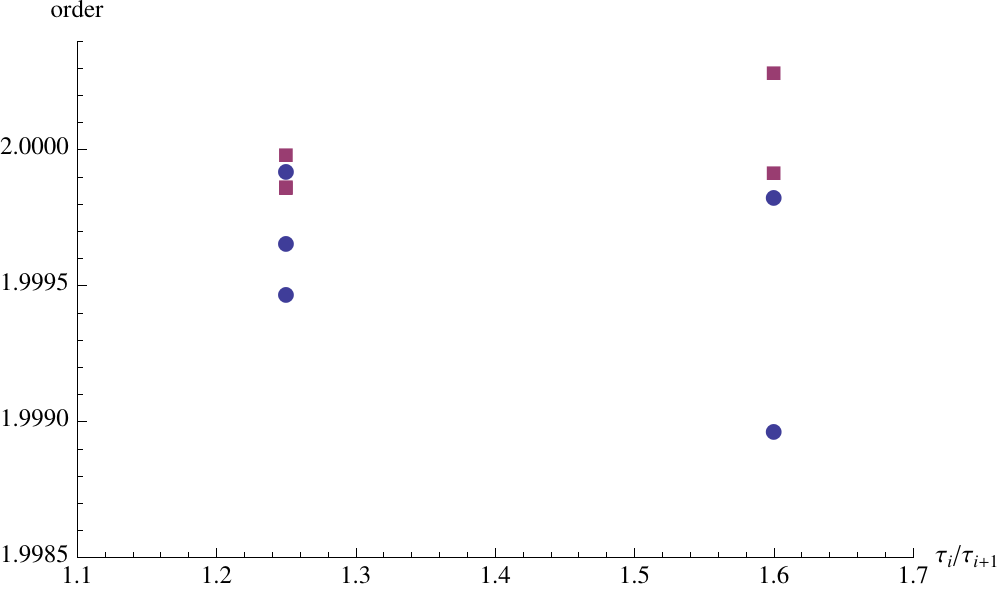}
\caption{Order approximation of MS type splittings as the function of $\tau_i/\tau_{i+1}.$
}\label{regularfig1}
\end{figure}
\end{center}

\newpage

\section{Discussion and perspectives}
We presented symbolic calculations for orders of PDE solving methods. Our motivation is to predict the order in the case when beside numerical procedures of certain order operator splitting is also used. We calculated the order of combined methods applied for nonlinear PDE-s like \eqref{uder}, where a bounded linear operator and a nonlinear operator is present. We presented numerical calculations on a test problem with the diffusion operator which is an unbounded linear operator. Although the results are in accordance with our theoretical results the methods used in section 5 are strictly correct for the case of bounded linear operator. The results of section 6 indicates that the combined method inherits the smaller one of the  order of the splitting and the numerical method, the extension of the methods used in section 5 to unbounded operators is not obvious. Our main focus is on reaction-diffusion equations so we plan to find a method which allows us to repeat the results of this paper for unbounded linear operators. 

We also intend to extend the results of \cite{csomosfarago} to nonlinear problems. 
Besides we are working on reaction-diffusion simulations on the sphere. In the future we are going to apply the presented methods to practical problems e. g. in the simulation of combustion.

\section{Acknowledgements}
This work is connected to the scientific program of the "Development
of quality-oriented and harmonized R+D+I strategy and functional model
at BME" project. This project is supported by the New Hungary
Development Plan (Project ID: TÁMOP-4.2.1/B-09/1/KMR-2010-0002). 
The present research has partially been supported by the National Science Foundation, Hungary (No. K84060).
I thank to Istv\'{a}n Farag\'{o} for the fruitful discussions.

\bibliographystyle{elsarticle-num}
\bibliography{<your-bib-database>}

\begin{thebibliography}{00}
\bibitem[Ablowitz and Zeppetella]{abl} Ablowitz, M. J.; Zeppetella, A.: Explicit solutions of Fisher's equation for a special wave speed, Bull. Math. Biol. 41, 835--840, 1979.

\bibitem[Csom\'{o}s and Farag\'{o}]{csomosfarago}
Csomós, P., Faragó, I.: Error analysis of the numerical solution of split differential equations, Math. Comp. Mod., 48, Issue 7--8, 1090--1106, 2008.

\bibitem[\'{E}rdi and T\'{o}th]{erditoth} \'{E}rdi, P., T\'{o}th, J., Ch. 6: Mathematical Models of Chemical Reactions, Princeton University Press, Princeton, N.J. 1989.

\bibitem[Farag\'{o} and Havasi]{faragohavasi} Farag\'{o}, I., Havasi \'{A}.: The mathematical background of operator splitting and the effect of non-commutativity, LEcture Notes in Computer Science, Vol. 2179/2001, 264--27, 2001.

\bibitem[Fisher]{fisher} Fisher, R. A.: The genetical theory of natural selection. Oxford University Press, Oxford, 1930.

\bibitem[Hundsdorfer and Verwer]{verwer} Hundsdorfer, W., Verwer, J. G.: Numerical solution of the time-independent advection-diffusion-reaction equations, Springer-Verlag, Berlin, 2003.

\bibitem[Ladics]{ladics} Ladics, T.:  Analysis of the splitting error for advection-reaction problems in air pollution models, Id\H{o}j\'{a}r\'{a}s, Quarterly Journal of the Hungarian Meteorological Service Vol. 109, No. 3, 173--188, July--September 2005.

\bibitem[Lagzi et al.]{lagzi1} Lagzi, I., Tomlin, A. S., Tur\'{a}nyi, T., Haszpra, L., M\'{e}sz\'{a}ros, R., Berzins, M.: Modeling Photochemical Air Pollution in Hungary Using an Adaptive Grid Model, 'Air Pollution Modelling and Simulation', pp. 264-273, editor: B. Sportisse, Springer, Berlin, 2002, ISBN 3-540-42515-2.
  
\bibitem[Lagzi et al.]{lagzi2} Lagzi, I., Tomlin, A. S., Tur\'{a}nyi, T., Haszpra, L., M\'{e}sz\'{a}ros, R., Berzins, M.: The Simulation of Photochemical Smog Episodes in Hungary and Central Europe Using Adaptive Gridding Models. Lecture Notes in Computer Science, 2074, 67-76, 2001.

\bibitem[Lagzi et al.]{lagzi3} Lagzi, I., Tomlin, A. S., Tur\'{a}nyi, T., Haszpra, L.: Photochemical air pollutant formation in Hungary using an adaptive gridding technique, Int. J. Environment and Pollution, 36, 44-58, 2009.

\bibitem[Lanser and Verwer]{lanserverwer} Lanser, D., Verwer, J. G.: Analysis of operator splitting for advection--diffusion--reaction problems from air pollution modeling, Journal of Computational and Applied Mathematics, Vol. 111 Issue 1--2, 201--216, Nov. 15 1999.

\bibitem[Larsson and Thom\'{e}e]{thomee} Larsson, S., Thom\'{e}e, V.: Partial Differential Equations with Numerical Methods, Springer--Verlag Berlin Heidelberg 2003.

\bibitem[Komornik]{komornik} Komornik, V.: Pr\'{e}cis d'analyse r\'{e}elle, tome I -- published by Ellipses -- copyright \'{E}dition Marketing S. A. 2001. Hungarian edition Typotex, Budapest 2003.

\bibitem[Marchuk]{marchuk}
Marchuk, G. I.: Methods of splitting, Nauka, Moscow, 1988. (in Russian)

\bibitem[Murray]{murray}:
Mathematical Biology, Ch. 11. 3rd edition in 2 volumes: Mathematical Biology: I. An Introduction (551 pages) 2002; Mathematical Biology: II. Spatial Models and Biomedical Applications (811 pages) 2003.

\bibitem[Sanz-Serna]{sanz} Sanz-Serna, J. M.:The State of the Art in Numerical Analysis, chapter Geometric Integration,
121--143, Clarendon Press, Oxford, 1997.

\bibitem[Simpson and Landman]{simpson} Simpson, M. J., Landman, K. A.: Characterizing and minimizing the split error for Fisher's equation, Appl. Math. Lett., 19 (7), 604--612, 2006.

\bibitem[Sportisse]{sportisse} Sportisse, B.: An analysis of operator splitting techniques in the stiff case. J. Comp. Phys. 161 (1), 140--168, 2000.

\bibitem[Strang]{strang}
Strang, G.: On the construction of different splitting schemes, SIAM J. Numer. Anal. 5 (3) 506--517, 1968.
\end{thebibliography}

\end{document}